\documentclass{article}

\usepackage[pdftex]{graphicx}
\usepackage{braket}
\usepackage{amsmath,amssymb}
\usepackage{mathrsfs}
\usepackage[top=30truemm,bottom=30truemm,left=40truemm,right=40truemm]{geometry}
\usepackage{indentfirst}

\usepackage{bm}
 \usepackage{amssymb}
 \usepackage{amsbsy}
 \usepackage{amscd}
 \usepackage{amsmath}
 \usepackage{amsthm}
\usepackage{color}

\newcommand{\EE}{\mathbb E}
\newcommand{\VV}{\mathbb V}

\begin{document}
\title{Asymptotic Behavior of Free Energy When Optimal Probability Distribution Is Not Unique}
\author{Shuya Nagayasu and Sumio Watanabe
\\
Department of Mathematical and Computing Science
\\
Tokyo Institute of Technology, Japan}
\date{}
\maketitle

\begin{abstract}
Bayesian inference is a widely used statistical method. The free energy and generalization loss,
which are used to estimate the accuracy of Bayesian inference, are known to be small 
in singular models that do not have a unique optimal parameter.
However, their characteristics are not yet known when there are multiple optimal probability distributions. 
In this paper, we theoretically derive the asymptotic behaviors of the generalization loss and free energy in the case that the optimal probability distributions are not unique and show that they contain asymptotically different terms from those of the conventional asymptotic analysis. 
\end{abstract} 

\section{Introduction}

In statistical learning theory, a probability distribution which generates a sample is called 
a true distribution and one with a parameter is called a statistical model
or a learning machine. An probability distribution is estimated by applying a
training algorithm to a statistical model. 
Then, the difference between the true distribution and the estimated one is defined by some measure, for example, the Kullback-Leibler (KL) divergence. 
In practical applications, the true distribution is unknown, hence the free energy and the generalization loss, which give the relative difference of KL divergence, are used to evaluate the estimated one.

The theoretical values of the free energy and the generalization loss strongly depend on
the geometrical situations of the true distribution and a statistical model. 
A statistical model is called regular if the parameter which minimizes the KL divergence of 
a true distribution and the statistical model is unique and Hessian matrix of the KL divergence 
at the minimum point is regular. For the regular case, the asymptotic behavior of the generalization loss was revealed by Akaike\cite{Akaike1974}, while that of the free energy was revealed by Schwarz\cite{Schwarz}. These results have been applied to statistical model selection criteria, i.e. Akaike(AIC), Bayesian(BIC), Deviance(DIC)\cite{Spiegelhalter}, and Adjusted Bayesian(ABIC)\cite{Akaike1980}. 

If a statistical model is not regular, then it is called singular. There are two different singular cases. One is that multiple optimal parameters exist, but the probability distributions of the optimal parameters are the same and unique. This case occurs when a statistical model such as a neural network or a normal mixture is redundant for to the true distribution. From here on, we will refer to this case as the multiple parameters and unique distribution case. 

The other case is when there are multiple optimal parameters that give different multiple probability distributions. 
This case occurs when there exist different optimal statistical models have the same KL divergence from the true distribution. We will call this case the multiple parameters and multiple distribution case. 

In the former singular case, the asymptotic behavior of the free energy was revealed by Watanabe\cite{S Watanabe1999}, where it was applied to the statistical model criteria WAIC\cite{S Watanabe2010a} and 
WBIC\cite{S Watanabe2013}. In contrast, asymptotic behavior in the latter case remains unknown. There, we decide to study the multiple parameters and multiple distribution case, and clarify the asymptotic behavior of the free energy. In particular we find that there is a new term that is not included in the previous statistical learning theory. 

\section{Main Result}

In this section, we explain the framework of Bayesian inference and show the main result of this paper. 

We assume that, from a true distribution $q(x)$, a set of $n$ independent random variables 
$X^n = (X_1, X_2, \cdots, X_n)$ and $X$ are generated, 
where $X_i \in \mathbb{R}^N (0 \leq i \leq N)$ and $X^n$ and $X$ are
independent.
Let $\EE_X[\;\;]$ and $\EE[\;\;]$ denote expectation values 
over $X$ and $X^n$. 
A statistical model and a prior are represented by probability density functions 
$p(x|w)$ and $\varphi(w)$, 
$ w \in W \subset R^d$, respectively. 
In Bayesian inference, for a given inverse temperature $\beta>0$, 
the marginal likelihood $Z_n(\beta)$, 
the posterior distribution $p(w|X^n)$, 
and the predictive distribution $p(x|X^n)$ are defined as
\begin{align*}
Z_n(\beta) &= \int \prod_{j = 1}^{n}p(X_j|w)^{\beta} \varphi(w) dw, \tag{1}\\ 
p(w|X^n) &= \frac{1}{Z_n} \prod_{j = 1}^{n}p(X_j|w)^{\beta} \varphi(w), \tag{2}\\
p(x|X^n) &= \int p(x|w) p(w|X^n) dw. \tag{3}
\end{align*}
Note that $\beta=1$ results in the conventional Bayesian estimation. 
The free energy $F_n(\beta)$ and the generalization loss $G_n(\beta)$ are 
defined as 
\begin{align*}
F_n(\beta) &= -\frac{1}{\beta}\log Z_n(\beta), \tag{4}
\\
G_n(\beta) &= - \EE_X[\log p(X|X^n)].   \tag{5}
\end{align*}
The set of optimal parameters $W_0$ is defined as the set of all parameters that minimize the KL divergence of $q(x)$ and $p(x|w)$, 
\begin{align*}
W_0=\{w\in W\;;\;
\int q(x) \log \frac{q(x)}{p(x|w)} dx \mbox{ is minimized. }.  \}\tag{6}
\end{align*}
The log density ratio function for $w_0\in W_0$ and $w\in W$ is defined as
\begin{align*}
f(x,w_0,w) = \log \frac{p(x|w_0)}{p(x|w)}. \tag{7}
\end{align*}
Then, by defining 
\[
L(w)=-\EE_X[\log p(X|w)],
\]
it follows that
\[
\EE_X[f(X,w_0,w)]=L(w)-L(w_0). 
\]
It is said that the log density ratio function has a relatively finite variance if the following condition is satisfied, 
\begin{align*}
\exists c_0>0, \forall w\in W, 
 \forall w_0 \in W_0, \; \;\EE_X[f(X,w_0,w)] \geq c_0 \EE_X[f(X,w_0,w)^2], \tag{8}
\end{align*}
It was proved in \cite{S Watanabe1999} that, if the log density ratio function
has a relatively finite variance, 
then the free energy and the generalization loss have asymptotic expansions, 
\begin{align*}
\EE[F_n(1)] &= nL(w_0) + \lambda \log n - (m - 1) \log \log n + O(1), \tag{9} \\
\EE[G_n(1)] &= L(w_0) + \frac{\lambda}{n} - \frac{m - 1}{n \log n} + o\left(\frac{1}{n\log n}\right), \tag{10}
\end{align*}
even if the Fisher information matrix is singular, 
where $\lambda > 0$ is a rational number called the real log canonical threshold (RLCT), 
and $m \geq 1$ is a natural number called the multiplicity. 
It was revealed that the assumption (8) is satisfied if a parameter $w_0$ exists that satisfies $p(x|w_0) = q(x)$, 
or if the Fisher information matrix is positive definite at a unique $w_0$.

However, in general, the assumption (8) is not always satisfied. For example, it is known that, if 
the optimal probability distribution is not unique, that is to say,
\begin{align*}
\exists w_{01}, w_{02} \in W_0, \; \; p(x|w_{01}) \neq p(x|w_{02}), \tag{11}
\end{align*}
then assumption (8) does not hold (See Lemma 1). Hence neither eq. (10) nor eq. (11) does
not holds in general. Previous research \cite{S Watanabe2010b} treated the case the assumption (8) does not hold, but the optimal probability distribution is unique. 

In this paper, we study the case in which the optimal probability distribution is not unique, as shown in eq.(11), and show that the asymptotic behaviors of the free energy and the generalization loss are given by 
\begin{align*}
\EE[F_n(1)] &= nL(w_0) - \mu \sqrt{n} + \hat{\lambda} \log n - (\hat{m} - 1) \log \log n + O(1), \tag{13} \\
\EE[G_n(1)] &= L(w_0) - \frac{\mu}{2\sqrt{n}} + o\left(\frac{1}{\sqrt{n}}\right), \tag{14}
\end{align*}
where $\mu, \hat{\lambda}, \hat{m}$ are real numbers satisfying $\mu \geq 0$, $\hat{\lambda} > 0$, and $ \hat{m} \geq 1$. 
Since $\mu>0$, 
these results show that both the free energy and the generalization loss are made smaller 
in Bayesian estimation if the optimal probability distribution is not unique. 

\section{Definitions and Notation}

In this section, we summarize several definitions and notation. 
We estimate a probability distribution $q(x)$ by Bayesian inference using a statistical
model $p(x|w)$ and prior distribution $\varphi(w)$. It is assumed that 
the set of parameters $W$ is compact and that $p(x|w)$ is continuous for $w$.
$W_0 \subset W$ is set of all optimal parameters that minimizes $L(w)$. 
In this case, the log density ratio function does not have a relatively finite variance. 
In fact, if eq. (8) holds,
 there exists a positive real number $c_0>0$ such that 
for any $w_{01},w_{02}\in W_0$, 
\begin{align*}
0 = L(w_{01}) - L(w_{02}) = \EE_X[f(X,w_{01},w_{02})] \geq c_0 E_X[f(X,w_{01},w_{02})^2], 
\end{align*}
and it follows that $f(x,w_{01},w_{02}) = 0$ for all $x$,  which means $p(x|w_{01}) = p(x|w_{02})$.
\vskip5mm
\noindent{\bf Problem treated in this paper}. 
We study the case, 
\begin{align*}
\exists w_{01}, w_{02} \in W_0 \; \; p(x|w_{01}) \neq p(x|w_{02}), 
\end{align*}
resulting that the assumption, eq.(8), does not hold. We assume that 
the set $W_0$ can represented by a disjoint union of $W_{0i} (i \in I)$ which satisfies
\begin{align*}
\cup_{i \in I} W_{0i} &= W_0, \\
W_{0i} \cap W_{0j} &= \varnothing (i \neq j),
\end{align*}
and in each subset $W_{0i}$, the assumption holds. 
Hence, for an arbitrary $w_{0i}\in W_{0i}$, $p(x|w_{0i})$ is the 
same probability distribution. 
\vskip5mm
The empirical log loss function is defined by 
\begin{align*}
L_n(w) &= -\frac{1}{n} \sum_{j = 1}^{n} \log p(X_j|w).\tag{15} 
\end{align*}
Accordingly, $L(w_{0i})$  does not depend on $i\in I$, and the random variable depends on $i$:
\[
L_{ni} = L_{n}(w_{0i}).
\]
The log loss function $L(w)$ and the average error function $K(w)$ are defined by
\begin{align*}
L(w) &= -\EE_X[\log p(X|w)],\tag{16}\\
K(w) &= \EE_X[f(X,w_{0i},w)]=L(w)-L(w_{0i}). \tag{17}
\end{align*}
Accordingly $K(w)$ does not depend on $w_{0i}$
and K(w) satisfies
\begin{align*}
K(w_{0i}) = 0, K(w) \geq 0. \tag{18}
\end{align*}
The empirical error function $K_{ni}(w)$ is defined by 
\begin{align*}
K_{ni}(w) = \frac{1}{n} \sum_{j = 1}^{n}f(X_j,w_{0i},w).\tag{19}
\end{align*}
It follows that
\[
L_n(w)=K_{ni}(w)+L_{ni}. 
\]
\newpage

\section{Asymptotic Properties of Free Energy }

In this section, we discuss the case in which the number of optimal probability distributions is finite. If the number of optimal probability distributions is a finite natural number m, 
the elements of the index set $I$ are natural numbers $(1,2,\cdots m)$. We will show that
\begin{align*}
\EE[F_n(\beta)] = nL(w_0) - \sqrt{n}\mu + 
\frac{1}{\beta}
(\hat{\lambda}\log n - (\hat{m} - 1)\log \log n) + O(1), \tag{20}
\end{align*}
holds. If $\EE[G_n(1)]$ has an asymptotic expansion, then
\begin{align*}
\EE[G_n(1)] = L(w_0) -\frac{\mu}{2\sqrt{n}}+ o\left(\frac{1}{\sqrt{n}}\right) \tag{21}
\end{align*}
holds.
\vskip5mm
\noindent{\bf Example}. The following is an example for the case that 
the optimal probability distribution is not unique. We suppose supervised learning of 
$q(y|x)$ using a statistical model $p(y|x,a,b)$. The variables $a,b$ are parameters. 
We suppose true distribution $q(y|x), q(x)$ is
\begin{align*}
q(y|x) &= \frac{1}{\sqrt{2 \pi}}\exp\left(-\frac{(y - f(x))^2}{2}\right), \\
q(x) &= \left\{\begin{array}{cc}
    \frac{1}{4} & (-2 \leq x \leq 2) \\
    0 & (\mbox{otherwise})
\end{array}\right.,
\\
f(x) &=\left\{ \begin{array}{cc}
    x + 2 & (-2 \leq x < -1) \\
    1 & (-1 \leq x < 1)\\
    -x + 2 &  (1 \leq x \leq 2)
\end{array}\right..
\end{align*}
We also suppose a statistical model $p(y|x,a,b)$,
\begin{align*}
p(y|x,a,b) &= \frac{1}{\sqrt{2 \pi}}\exp\left(-\frac{(y - \sigma(ax + b))^2}{2}\right) \\
\sigma(x) &= \frac{1}{1 + \exp(-x)}.
\end{align*}
In this situation, the KL divergence between $q(y|x)$ and $p(y|x,a,b)$ can be calculated as 
\begin{align*}
KL(q(y,x)|p(y,x|a,b)) &= \int q(y,x) \log \frac{q(y|x)}{p(y|x,a,b)}
= \int q(y|x)q(x) \log \frac{q(y|x)}{p(y|x,a,b)}, \\
&=\frac{1}{2}  \int q(y|x)q(x) (f(x) - \sigma(ax + b))^2.
\end{align*}
Note that 
$q(x)$ and $f(x)$ are even functions. In addition, $\sigma(ax + b)$ and $\sigma(-ax + b)$ are line symmetric on $x = 0$. Therefore, $KL(q(y|x)|p(y|x,a,b)) = KL(q(y|x)|p(y|x,-a,b))$ holds. Thus, 
 we can find two optimal parameters as $(a_0,b_0)$ and $(-a_0,b_0)$. These two points satisfy $p(y,x|a_0,b_0) \neq p(y,x|-a_0,b_0)$. A numerical calculation shows that the optimal parameters in this case are
$ (5.13 , 7.71)$, $(-5.13 , 7.71)$.
\vskip5mm
To show eq.(20), let us represent the set of parameters as 
\[
W=\cup_{i=1}^m W_i,
\]
where $W_{0i}\subset W_i$. 
The marginal likelihood of each domain $Z_n^{(i)}$ is also defined by
\begin{align*}
Z_n^{(i)} = \int_{w \in W_i} \exp(-n\beta L_n(w))\varphi(w)dw. \tag{22}
\end{align*}
By using $L_n(w)=K_{ni}(w)+L_{ni}$, we have
\begin{align*}
Z_n^{(i)} = \exp(-n\beta L_{ni})\int_{w \in W_i} \exp(-n\beta K_{ni}(w))\varphi(w)dw. \tag{23}
\end{align*}
The normalized marginal likelihood of each domain $Z_n^{0i}$ is defined by
\begin{align*}
Z_n^{0i} = \int_{w \in W_i} \exp(-n\beta K_{ni}(w))\varphi(w)dw. \tag{24}
\end{align*}
In accordance with these definitions, the marginal likelihood is given by 
\begin{align*}
Z_n = \int_{w \in W} \exp(-n\beta L_n(w))\varphi(w)dw = \sum_{i = 1}^{m}
\exp(-n\beta L_{ni})Z_n^{0i}. \tag{25}
\end{align*}
The normalized marginal likelihood $Z_n^{0i}$ can be divided into $Z_n^{1i}$ and $Z_n^{2i}$  as 
\begin{align*}
&Z_n^{0i} = Z_n^{1i} + Z_n^{2i},\tag{22}\\
&Z_n^{1i} = \int_{w \in W_i,K(w) < \epsilon} \exp(-nK_{ni}(w))\varphi(w)dw, \tag{26}\\
&Z_n^{2i} = \int_{w \in W_i,K(w) \geq \epsilon} \exp(-nK_{ni}(w))\varphi(w)dw. \tag{27}
\end{align*}
By taking $\epsilon>0$ to be a monotonically decreasing function of n which satisfies
\begin{align*}
\lim_{n \rightarrow \infty} \epsilon = 0,\\
\lim_{n \rightarrow \infty} \sqrt{n} \epsilon = \infty,
\end{align*}
it can be proved that
$Z_n^{2i}=O_p(\exp(-\sqrt{n}))$\cite{S Watanabe1999}.
The optimal parameter set in each $W_i$ is $W_{0i}$;
therefore, in this set, the optimal probability distribution is only $p(x|w_{0i})$. By assumption, 
\begin{align*}
\exists c_0 \forall w \in W_i \forall w_{0i} \; \; E_X[f(x,w_{0i},w_i)] \geq c_0 E_X[f(x,w_{0i},w_i)^2].
\end{align*}
Then, applying the results of \cite{S Watanabe1999}, by using the resolution theorem on each $W_i$, 
there exists a measure $du_i^*$ such that 
\begin{align*}
Z_n^{(1i)} = \frac{(\log n)^{m_i-1}}{n^{\lambda_i}}\int du^*_i\int t^{\lambda_i - 1} \exp(-\beta t + \beta\sqrt{t}\xi_{ni}(u_i))dt + o_p\left(\frac{(\log n)^{m_i-1}}{n^{\lambda_i}}\right). \tag{28}
\end{align*}
In this equation, $\lambda_i$ and $m_i$ are respectively the real log canonical threshold and
the multiplicity of zeta functions, and $\xi_{ni}(u_i)$ is an empirical process
that converges in distribution to a Gaussian process. 

Thus, the marginal likelihood is
\begin{align*}
Z_n &= \sum_{i = 1}^{m}\exp(-n\beta L_{ni})Z_n^{0i}\\
      &= \sum_{i = 1}^{m}\exp(-n\beta L_{ni})
\frac{(\log n)^{m_i-1}}{n^{\lambda_i}}\left(
\int du^*_i\int t^{\lambda_i - 1} 
\exp(-\beta t + \beta\sqrt{t}\xi_{ni}(u_i))dt + o_p\left(1\right)\right). \tag{29}
\end{align*}
From (4), the free energy is given by   
\begin{align*}
F_n(\beta) &= -\frac{1}{\beta}\log Z_n(\beta) 
\\
&=f_1+f_2+f_3+ o_p(1),\tag{30}
\end{align*}
where
\begin{align*}
f_1&= -\frac{1}{\beta}\log\left(\sum_{k = 1}^{m}e^{(-n\beta L_{ni})}\right),
\\
f_2&=- \frac{1}{\beta}\log\left(\frac{\sum_{i = 1}^{m}e^{(-n\beta L_{ni} - \lambda_i \log n + (m_i - 1)\log \log n)}}{\sum_{k = 1}^{m}e^{(-n\beta L_{ni})}}\right), \\
f_3&=- \frac{1}{\beta}\log\left(\frac{\sum_{i = 1}^{m}
e^{(-\Theta(\beta,\xi_{ni})-n\beta L_{ni} - \lambda_i \log n + (m_i - 1)\log \log n)}}
{\sum_{i = 1}^{m}e^{(-n\beta L_{ni} - \lambda_i \log n + (m_i - 1)\log \log n)}}\right).
\end{align*}
In the above equations equation, we have used the notation,
\begin{align*}
&\Theta(\beta,\xi_{ni})\\ &\qquad = - \log \left(\int du_i^{*} \int_{0}^{\infty} dt \, t^{\lambda_i - 1} \exp(- \beta t + \beta\sqrt{t}\xi_{ni}(u_i))\right).\tag{31}
\end{align*}
In the following, we examine the asymptotic behaviors of the three terms eq.(30). 
First, to study $f_1$, we define $i_{max}$ and $Y$ by 
\begin{align*}
i_{max} &= \underset{i}{\mbox{argmax}} \; (- L_{ni}),\\
Y &= -n\beta L_{ni_{max}}.
\end{align*}
From the definition,
\begin{align*}
\log\left(\sum_{i = 1}^{m}e^{(-n\beta L_{ni})}\right) - Y = \log\left(1 + \sum_{i \neq i_{max}}^{m}e^{(-n\beta L_{ni} - Y)}\right).\tag{32}
\end{align*}
Since $Y+n\beta L_{ni}\geq 0$, 
\begin{align*}
0 < \log\left(\sum_{i = 1}^{m}e^{(-n\beta L_{ni})}\right) - Y \leq m\log 2. 
\end{align*}
Therefore,
\[
f_1=Y/\beta+O_p(1)= -n L_{ni_{max}}+O_p(1).
\]
Note that
the average of $L_{ni}$ is $L(w_{0i})$, which does not depend on $i$. 
Let us define a random variable, 
\begin{align*}
\mathscr{L}_{n}(w_{0i})\equiv \sqrt{n} (-L_{n}(w_{0i})+L(w_{0i})). \tag{33}
\end{align*}
By using the central limit theorem, 
$  \mathscr{L}_{n}(w_{0i})$ $(i=1,2,...,m)$ 
converges  in distribution to an m-dimensional Gaussian random variable  
$\mathscr{L}(w_{0i})$ on $w_{0i}\in W_0$ 
whose average is zero and 
variance-covariance matrix  $V=(V_{ij})$ is
\begin{align*}
V_{ij} = \EE_X[(\log p(X|w_{0i}) + L(w_0))(\log p(X|w_{0j}) + L(w_0))]. \tag{34}
\end{align*}
Then we obtain
\[
f_1=n L(w_0) -
\sqrt{n}\underset{w_0 \in W_0}{\mbox{max}} \mathscr{L}_n(w_0)+O_p(1). \tag{35}
\]
Using $\mathscr{L}(w_0)$, the asymptotic behavior of its average is given by
\[
\EE[f_1]=n L(w_0) -
\EE[\sqrt{n} \underset{w_0 \in W_0}{\mbox{max}} \mathscr{L}(w_0)]+O(1). \tag{36}
\]
Now, let us examine the second term $f_2$ in eq.(30). 
If there exist multiple $i_{max}$, we define $i_{max}$ as the $i$ whose RLCT is smallest, 
and if  the RLCTs are the same, we define $i_{max}$ as the $i$ whose multiplicity is biggest. 
Using $i_{max}$ so defined, the asymptotic behavior of the second term $f_2$ 
is given by 
\begin{align*}
-\beta f_2
&=\log\left(\frac{\sum_{i = 1}^{m}e^{(-n\beta L_{ni} - \lambda_i \log n + (m_i - 1)\log \log n)}}{\sum_{i = 1}^{m}e^{(-n\beta L_{ni})})}\right)\\
&\qquad = - \lambda_{i_{max}} \log n + (m_{k_{max}} - 1)\log \log n 
\\
&\qquad 
+ \log\left(\frac{1 + \sum_{k \neq i_{max}}^{m}
e^{(-n\beta L_{nk} - Y - (\Delta\lambda_k) \log n + (\Delta m_k)\log \log n)}}
{1 + \sum_{k \neq i_{max}}^{m}e^{(-n\beta L_{ni} - Y)})}\right)\\
&\qquad = - \lambda_{i_{max}} \log n + (m_{i_{max}} - 1)\log \log n 
+ \log\left(\frac{1 + o_p(e^{-n\beta a})}{1 + o_p(e^{-n\beta a})}\right)\\
&\qquad = - \lambda_{i_{max}} \log n + (m_{i_{max}} - 1)\log \log n 
+ o_p(e^{-n\beta a}), \tag{37}
\end{align*}
where $\Delta\lambda_i$ and $ \Delta m_i$ are $\lambda_i - \lambda_{i_{max}}$
 and $m_i - m_{i_{max}}$ respectively, and $a\geq0 $ is the difference between $Y$ and the 
second biggest $-nL_{ni}$.
\begin{align*}
\EE[f_2]= \frac{1}{\beta}\left(
\lambda_{i_{max}} \log n - (m_{i_{max}} - 1)\log \log n 
\right)
+ O(1). \tag{38}
\end{align*}
Next, let us study study the third term $f_3$ in eq.(30). 
We define a random variable $a_i$ as follows.
\begin{align*}
a_i = \frac{e^{-n\beta L_{ni} - \lambda_i \log n + (m_i - 1)\log \log n}}{\sum_{i = 1}^{m}e^{(-n\beta L_{ni} - \lambda_i \log n + (m_i - 1)\log \log n)}}. \tag{39}
\end{align*}
The sum of $a_i$ over $i$ is 1. Using $a_i$ is 1. We can describe the third term using $a_i$, so we can describe the free energy as follows.
\begin{align*}
F_n(\beta) = nL(w_0) &- \sqrt{n}\underset{w_0 \in W_0}{\mbox{max}} \mathscr{L}_n(w_0) + \frac{\lambda_{i_{max}}}{\beta} \log n- \frac{(m_{i_{max}} - 1)}{\beta}\log \log n\\ &- \log \left(\sum_{i = 1}^{m}a_i e^{- \Theta(\beta,\xi_{ni})} \right) + O_p(1). \tag{40}
\end{align*}
\\
Lastly, in order to derive the asymptotic behavior of
 the average $\EE[F_n]$, we show that the average of the random variable
\[
f_4\equiv \log \left(\sum_{i = 1}^{m}a_i e^{- \Theta(\beta,\xi_{ni})} \right) 
\]
is finite. By the Cauchy-Schwarz inequality,
\begin{align*}
-\frac{t + \mbox{sup}_u|\xi_{ni}(u_i)|^2}{2} \leq \sqrt{t}\xi_{ni}(u_i) \leq \frac{t + \mbox{sup}_u|\xi_{ni}(u_i)|^2}{2} \tag{41}
\end{align*}
holds. In the integral range of $u_i$,$[0,1]^d$, we have
\begin{align*}
&-\log \int du_i^* -\log \int dt \, t^{\lambda_i - 1}\exp(-\frac{\beta}{2} t) - \frac{\beta}{2}\underset{u_i \in [0,1]^d}{\sup}|\xi_{ni}(u_i)|^2 \\ &\leq \Theta(\beta,\xi_{ni}(u_i))\\ &\leq  -\log \int du_i^* -\log \int dt \, t^{\lambda_i - 1}\exp(-\frac{3\beta}{2} t) + \frac{\beta}{2}\underset{u_i \in [0,1]^d}{\sup}|\xi_{ni}(u_i)|^2. \tag{42}
\end{align*}
Using Jensen's inequality, 
\begin{align*}
f_4
&\qquad\geq \sum_{i = 1}^{m}a_i (- \Theta(\beta,\xi_{ni}(u)))\\ 
&\qquad\geq \sum_{i = 1}^{m}a_i\left(\log \int du_i^* + \log \int dt \, t^{\lambda_i - 1}\exp^{-\frac{3\beta}{2}t}- \frac{\beta}{2}\underset{u_i \in [0,1]^d}{\sup}|\xi_{ni}(u_i)|^2\right) \tag{43}
\end{align*}
holds. From the nature of emprical processes, $\xi_{ni}(u_i)$ converge in low to Gaussian processes 
$\xi_{i}(u_i)$ and
\begin{align*}
\lim_{n \rightarrow \infty}E[\underset{u_i \in [0,1]^d}{\sup}|\xi_{ni}(u_i)|^2] = \lim_{n \rightarrow \infty}E[\underset{u_i \in [0,1]^d}{\sup}|\xi_i(u_i)|^2]
\end{align*}
holds. In regard to the lower bound of $\log \left(\sum_{i = 1}^{m}a_i e^{- \Theta(\beta,\xi_{ni}(u))} \right)$, the only term that may diverge as a random variable is $\xi_{ni}(u_i)$, so $E[\log \left(\sum_{i = 1}^{m}a_i e^{- \Theta(\beta,\xi_{ni}(u))} \right)]$ can be shown to be bounded below. In addition, because $\underset{u_i \in [0,1]^d}{\sup}|\xi_{ni}(u_i)|^2 \geq 0$ holds, we have
\begin{align*}
f4 &\quad\leq \log \left(\max_{i}\int dt \, t^{\lambda_i - 1}e^{-\frac{\beta}{2}t}\int du_i^* \right)\left(m e^{\frac{\beta}{2}\sum_{i = 1}^{m} \mbox{sup}_u|\xi_{ni}(u_i)|^2}\right)\\ 
&\quad= \log \left (\max_{i}\int dt \, t^{\lambda_i - 1}e^{-\frac{\beta}{2}t}\int du_i^*\right) + \log m + \frac{\beta}{2}\sum_{i = 1}^{m} \mbox{sup}_u|\xi_{nk}(u_i)|^2. \tag{44}
\end{align*}
Hence as we did for the lower bound it can be shown that $\EE[\log \left(\sum_{i = 1}^{m}a_i e^{- \Theta(\beta,\xi_{ni}(u))} \right)]$ is bounded from above.
By summing up the above equations, the 
asymptotic behavior of $\EE[F_n]$ can be described as 
\begin{align*}
\EE[F_n(\beta)] &= nL(w_0) - \sqrt{n}\EE[\underset{w_0 \in W_0}{\mbox{max}} \mathscr{L}(w_0)]\\
 &+ \frac{1}{\beta}\sum_{k = 1}^{m}\alpha_i(\lambda_{i} \log n - (m_{i} - 1)\log \log n) + O(1),
 \tag{45}
\end{align*}
where $\alpha_i$ is the probability that $i=i_{max}$. By putting 
$\hat{\lambda}=\sum_{k = 1}^{m}\alpha_i\lambda_{i} $ and 
$\hat{m}=\sum_{k = 1}^{m}\alpha_i m_{i} $, we obtain
\begin{align*}
\EE[F_n(\beta)] &= nL(w_0) 
- \sqrt{n}\EE[\underset{w_0 \in W_0}{\mbox{max}} \mathscr{L}(w_0)]
+ \frac{\hat{\lambda}}{\beta}\log n - \frac{\hat{m} - 1}{\beta}\log \log n + O(1) \tag{46}
\end{align*}
holds.\\
When $\beta = 1$ holds, we have \cite{Amari,Levin}
\begin{align*}
\EE[G_n(1)] &= \EE[F_{n + 1}(1)] - \EE[F_n(1)]. \tag{47}
\end{align*}
Using this equation, and assuming that $\EE[G_n(1)]$ has an asymptotic 
expansion, we find that 
\begin{align*}
\EE[G_n(1)] =  L(w_0) -\frac{1}{2\sqrt{n}}\EE[\underset{w_0 \in W_0}{\mbox{max}} \mathscr{L}(w_0)] + o\left(\frac{1}{\sqrt{n}}\right).  \tag{48}
\end{align*}

\section{Experiment}
In this section, we show the results of an experiment for the case when the optimal probability distribution is not unique.We set the true distribution as
\begin{align*}
&q(y|x) = \frac{1}{\sqrt{2 \pi (0.2)^2}}\exp\left(-\frac{(y - f(x))^2}{2(0.2)^2}\right). \\
&f(x) = \begin{cases}
    x + 2 & (-2 \leq x < -1) \\
    1 & (-1 \leq x < 1)\\
    -x + 2 &  (1 \leq x \leq 2)
\end{cases}
\\
&q(x) = \begin{cases}
    \frac{1}{4} & (-2 \leq x \leq 2) \\
    0 & \mbox{otherwise}.
\end{cases}
\end{align*}
We use following statistical model and prior distributions. 
\begin{align*}
&p(y|x,a,b) = \frac{1}{\sqrt{2 \pi (0.2)^2}}\exp\left(-\frac{(y - \sigma(ax + b))^2}{2(0.2)^2}\right). \\
&\sigma(x) = \frac{1}{1 + \exp(-x)}.
\\
&\varphi(a) = \begin{cases}
    \frac{1}{20} & (0 \leq x \leq 20) \\
    0 & \mbox{otherwise}
\end{cases}
\\
&\varphi(b) = \begin{cases}
    \frac{1}{40} & (-20 \leq x \leq 20) \\
    0 & \mbox{otherwise}.
\end{cases}
\end{align*}
This statistical model has two optimal parameters 
\begin{align*}
w_{01} = (5.13 , 7.71), w_{02} = (-5.13 , 7.71).
\end{align*}
At these points,
\begin{align*}
p(x|w_{01}) \neq p(x|w_{02})
\end{align*}
holds. In this case, eq. (45) gives the theoretical asymptotic behavior of the free energy versus inverse temperature for $\beta = 1$. Note that the KL-divergence between q(y,x) and p(y,x|a, b) in each neighborhood of the optimal parameter is regular, so $\lambda = 1and m = 1$. The expectation of the maximum value of a 2-dimensional Gaussian distribution follows a 1-dimensional Gaussian distribution (see the appendix). We will show that the theoretical behavior of free energy obeys
\begin{align*}
\EE[F_n(1)] &= nL(w_0) 
- \sqrt{n}\sqrt{\frac{\VV[\log(p(y|x,a_0,b_0)) - \log(p(y|x,-a_0,b_0))]}{2\pi}}
+ \log n + O(1) \tag{49}.
\end{align*}
In eq. (49), $L(w_0)$ and the coefficient of $\sqrt{n}$ can be calculated by numerical integration. We used the average of $F_n$ calculated from the true distribution $q(y|x), q(x)$ as the experimental value of $\EE[F_n]$. The prior distribution $p(a)$, $p(b)$ does not have an effect on the asymptotic behavior. For this reason, we used equally spaced fixed values for integration. We compared this experimental values and theoretical values, except for the $O(1)$ term.

\newpage

\begin{figure}
        \centering
        \includegraphics[keepaspectratio, scale=0.3]{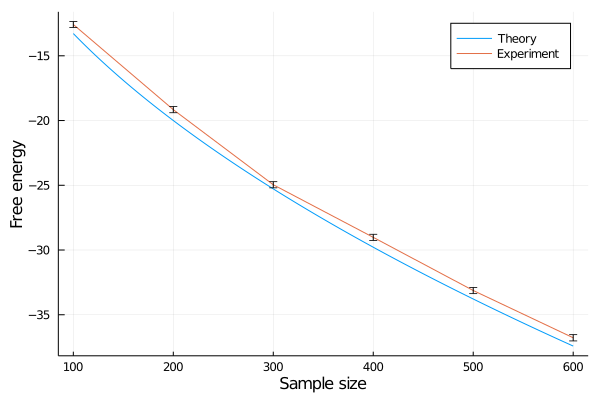}
        \caption{Experimental value and theoretical free energy depending on the sample size. The error bar is the SE of the average of free energy. The theoretical value of log likelihood$(nL(w_0))$ is subtracted from each value.}
        \label{figure1}
\end{figure}

\begin{figure}
        \centering
        \includegraphics[keepaspectratio, scale=0.3]{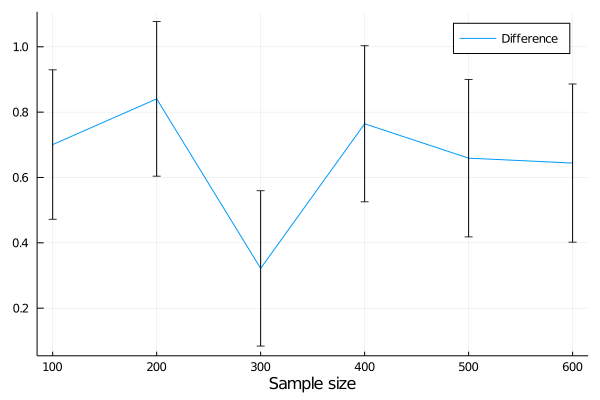}
        \caption{Difference between experimental value and theoretical value in Figure1. The error bar is the same as in Figure1.}
        \label{figure2}
\end{figure}

We calculated the experimental values of $\EE[F_n]$ whose sample sizes were
 $n=100$ to $600$ every $100$. The number of samples ranged from $10000$ to $60000$ in steps of $10000$ for each sample sizes. Figure1 compares the experimental and theoretical values. The experimental behavior of the free energy depending on the sample size is similar to the theoretical behavior. Figure2 shows the difference between the theoretical value and experimental value. This difference corresponds to $O(1)$ term.
This difference is remains on this order regardless of the sample size. The experimental results support the theoretical formula, eq.(45).

\section{Discussion}
We found that if the optimal probability distribution is not unique, the apparent bias or the variance gets smaller for a finite sample number $n$ corresponding to a Gaussian process determined by the log loss of  the optimal parameter set and the reduction converges to 0 asymptotically.This behavior can be explained qualitatively as follows: when there are two or more optimal probability distributions, the posterior distribution can be selected to be the nearest optimal probability distribution by bias of data, and this makes the generalization loss smaller than the average generation of data. As the number of samples and the bias increase, data generates averagely and generalization loss gets larger.\\
In this paper, we showed that the asymptotic behavior of the free energy and generalization loss are determined by $n^{\frac{1}{2}}$ and $n^{-\frac{1}{2}}$ order. Previous research\cite{S Watanabe2010b} provides concrete example in which there is a unique optimal probability distribution but assumption(8) does not hold. In that paper, the asymptotic behavior of the free energy and the generalization loss are determined by $n^{\frac{1}{3}}$ and $n^{-\frac{2}{3}}$ order, so we predict that the lowest order determining the asymptotic behavior of the free energy and generalization loss are $n^{\frac{1}{2}}$ and $n^{-\frac{1}{2}}$.\\
We showed that the asymptotic behaviors of the free energy and generalization loss are determined by the maximum value of a Gauusian process. The probability distribution of  the maximum value of a Gaussian process or multivariate normal distribution can not be calculated analytically, but an approximate calculation, called  the ``tube method''\cite{Kuriki} exists. There is also a a method for calculating the upper and lower bounds of the expectation of the maximum value of a Gauusian process, called ``chaining''\cite{Talagrand}.This maximum value is what determines the free energy and generalization loss in this paper.

\section{Conclusion}
We examined the case of when an important assumption in singular learning theory about the log density ration function is loosened. In this case there is a new term that is determined by a Gaussian process, whereby the generalization loss asymptotically increases as the size of the dataset increase. In the future, we should examine the asymptotic behavior of the generalization loss as a random variable, in particular the asymptotic equivalence of WAIC \cite{S Watanabe2010a} and WBIC \cite{S Watanabe2013} in this case, and in the case in which the assumption is completely removed.

\section*{Appendix}
We will derive the following equation.
\begin{align*}
\EE[{\mbox{max}} \mathscr{L}(w_0)] = \sqrt{\frac{\VV[\log(p(x|w_{01})) - \log(p(x|w_{02}))]}{2\pi}} \tag{A1}
\end{align*}
In this equation, $\mathscr{L}(w_0)$ is a 2-dimensional gaussian distribution which average is $0$ and variance-covariance matrix is  
\begin{align*}
V_{ij} = \EE[(\log p(x|w_{0i}) + L(w_0))(\log p(x|w_{0j}) + L(w_0))] \tag{A2}
\end{align*}
We define two random variables as 
\begin{align*}
z_1 = \mathscr{L}(w_{01}) - \mathscr{L}(w_{02}).\\
z_2 = \mathscr{L}(w_{01}) + \mathscr{L}(w_{02})
\end{align*}
The random variables $(z_1, z_2)$ are also from 2-dimensional Gaussian distribution whose average is $0$ and variance-covariance matrix is  
\begin{align*}
\left(
    \begin{array}{ccc}
      1 & -1  \\
      1 & 1 
    \end{array}
  \right)
\left(
    \begin{array}{ccc}
      V_{11} & V_{12}  \\
      V_{21} & V_{22} 
    \end{array}
  \right)&
\left(
    \begin{array}{ccc}
      1 & 1  \\
     -1 & 1 
    \end{array}
  \right)\\ &=
\left(
    \begin{array}{ccc}
      V_{11} + V_{22} - V_{12} - V_{21} & V_{11} - V_{22} + V_{12} - V_{21} \\
      V_{11} - V_{22} - V_{12} + V_{21} & V_{11} + V_{22} + V_{12} + V_{21}
    \end{array}
  \right).
\end{align*}
The marginal distribution about $z_1$ is a 1-dimensional Gaussian distribution whose average is $0$ and the variance is
\begin{align*}
V_{11} + V_{22} - V_{12} - V_{21}.
\end{align*}
According to eq.(16) and eq.(34), we have
\begin{align*}
V_{11} &+ V_{22} - V_{12} - V_{21}\\ 
&= \EE[(\log p(x|w_{01})(\log p(x|w_{01})] - L(w_0)^2 + \EE[(\log p(x|w_{02})(\log p(x|w_{02})]\\
& \qquad \qquad - L(w_0)^2 - 2(\EE[(\log p(x|w_{01})(\log p(x|w_{02})] - L(w_0)^2)\\
&= \EE[(\log p(x|w_{01}) - \log p(x|w_{02}))^2]\\
&= \VV[(\log p(x|w_{01}) - \log p(x|w_{02}))] \tag{A3}
\end{align*}
We define a random variable $z_3$ 
\begin{align*}
z_3 = \begin{cases}
    z_1 & z_1 \geq 0\\
    0 & z_1 < 0
\end{cases}
\end{align*}
By using $z_3$ we can describe the maximum value of $\mathscr{L}(w_0)$ in the following way,
\begin{align*}
{\mbox{max}} \mathscr{L}(w_0) = \mathscr{L}(w_{02}) + z_3. \tag{A4}
\end{align*}
Considering the average of $\mathscr{L}(w_{02})$ is $0$, we find that
\begin{align*}
\EE[{\mbox{max}}\mathscr{L}(w_0)] = \EE[z_3]. \tag{A5}
\end{align*}
$\EE[z_3]$ is the expectation of  a positive value in a Gaussian distribution. This integration of a Gaussian whose variance is $\sigma^2$ can be calculated as
\begin{align*}
\int_0^{\infty} \frac{x}{\sqrt{2 \pi}} \exp \left(-\frac{x^2}{2\sigma^2}\right) &= \left[-\frac{\sigma}{\sqrt{2 \pi}}\exp\left(-\frac{x^2}{2\sigma^2}\right)\right]_0^{\infty}\\
&= \frac{\sigma}{\sqrt{2 \pi}}. \tag{A6}
\end{align*}
From (A3), (A5), and(A6), we have 
\begin{align*}
\EE[{\mbox{max}}\mathscr{L}(w_0)] = \EE[z_3] = \sqrt{\frac{\VV[\log(p(x|w_{01})) - \log(p(x|w_{02}))]}{2\pi}}.
\end{align*}
Therefore (A1) holds.

\end{document}